\documentclass[a4paper,12pt,reqno,oneside]{amsart} 

\usepackage{setspace} 
\usepackage{typearea} 
\usepackage{amscd,amsmath,amssymb,verbatim} 

\usepackage{xr-hyper}
\usepackage{cite}
\usepackage[colorlinks,backref,final,bookmarksnumbered,bookmarks]{hyperref}
\usepackage[backrefs]{amsrefs}
\usepackage{graphicx}

\usepackage{ifthen}
\newcommand{\arxiv}[2][]{\ifthenelse{\equal{#1}{}}
{\href{http://arxiv.org/abs/#2}{\tt arXiv:#2}}
{\href{http://arxiv.org/abs/math/#2}{\tt arXiv:math.#1/#2}}}

\theoremstyle{plain}
\newtheorem{theorem}{Theorem}[section]
\newtheorem{mainthm}{Theorem}
 
\newtheorem{lemma}[theorem]{Lemma}
\newtheorem{proposition}[theorem]{Proposition}
\newtheorem{corollary}[theorem]{Corollary}
\newtheorem{maincor}[mainthm]{Corollary}

\theoremstyle{definition}
\newtheorem{remark}[theorem]{Remark}
\newtheorem*{remark*}{Remark}
\newtheorem{example}[theorem]{Example}

\def\x{\times}
\def\but{\setminus}
\def\emb{\hookrightarrow}

\def\phi{\varphi}

\renewcommand{\:}{\colon}

\def\R{\Bbb{R}}
\def\Z{\Bbb{Z}}
\def\xr#1{\xrightarrow{#1}} 
\makeatletter
\newcommand{\xR}[2][]{\ext@arrow 0359\Rightarrowfill@{#1}{#2}}
\newcommand{\xL}[2][]{\ext@arrow 0359\Leftarrowfill@{#1}{#2}}
\makeatother
\DeclareMathOperator{\id}{id}
\DeclareMathOperator{\const}{const}

\begin{document}
\title{Embeddability of joins and products of polyhedra}
\author{Sergey A. Melikhov}
\address{Steklov Mathematical Institute of Russian Academy of Sciences, Moscow, Russia}
\email{melikhov@mi-ras.ru}

\begin{abstract} We present a short proof of S. Parsa's theorem that there exists
a compact $n$-polyhedron $P$, $n\ge 2$, non-embeddable in $\R^{2n}$, such that $P*P$ 
embeds in $\R^{4n+2}$. 
This proof can serve as a showcase for the use of geometric cohomology.

We also show that a compact $n$-polyhedron $X$ embeds in $\R^m$, $m\ge\frac{3(n+1)}2$, if either
\begin{itemize}
\item $X*K$ embeds in $\R^{m+2k}$, where $K$ is the $(k-1)$-skeleton of the $2k$-simplex; or
\item $X*L$ embeds in $\R^{m+2k}$, where $L$ is the join of $k$ copies of the $3$-point set; or
\item $X$ is acyclic and $X\x\text{(triod)}^k$ embeds in $\R^{m+2k}$.
\end{itemize}
\end{abstract}
\maketitle

\section{Introduction}

It was shown by Flores, van Kampen and Gr\"unbaum \cite{Gr} that every $n$-dimensional join 
of $k_i$-skeleta of $(2k_i+2)$-simplexes does not embed in $\R^{2n}$ (see also 
\cite{M1}*{Examples 3.3, 3.5}, \cite{M2}, \cite{VS}).
Some other $k_i$-polyhedra with this property are constructed in \cite{M2}.

As noted by S. Parsa \cite{Pa}, it is implicit in a paper by Bestvina, Kapovich and Kleiner \cite{BKK} 
that if compact polyhedra $P^n$ and $Q^m$ both have non-zero mod $2$ van Kampen obstruction,
then $P*Q$ does not embed in $\R^{2(n+m+1)}$.
An $n$-dimensional polyhedron, non-embeddable in $\R^{2n}$ but with vanishing mod $2$ van Kampen 
obstruction was constructed by the author for each $n\ge 2$ \cite{M1}, settling a 1991 question by K. Sarkaria.
Using this example, Parsa proved the following, relying on extensive algebraic computations.

\begin{mainthm}[Parsa \cite{Pa}] \label{t1} For each $n\ge 2$ there exists a compact 
$n$-polyhedron $P$ such that $P$ does not embed in $\R^{2n}$ but $P*P$ embeds in $\R^{4n+2}$.
\end{mainthm}

In the present note we present a simple geometric proof of Parsa's theorem.

\begin{remark}
Parsa shows more generally that if compact polyhedra $P^n$ and $Q^m$ both have zero mod $2$ 
van Kampen obstruction, then $P*Q$ embeds in $\R^{2(n+m+1)}$ \cite{Pa}.
In fact, the proof of this more general assertion is implicit in our proof of Theorem \ref{t1}.

When precisely one of the mod $2$ van Kampen obstructions of $P$ and $Q$ is nonzero, 
our method can be used to check that $P*Q$ is non-embeddable in $\R^{2(n+m+1)}$ at least in
some cases (for instance, when $P=K$ and $Q=L$ in the notation of \S\ref{proof} below).
Parsa has a more precise result in this direction \cite{Pa}.
\end{remark}

By considering the cone $Q$ over the polyhedron $P$ in Theorem \ref{t1} we at once%
\footnote{See Corollary \ref{basic-a'} below and \cite{RS}*{Lemma 2.19 and Exercise 2.24(3)}.}
obtain

\begin{maincor} \label{t2} For each $n\ge 3$ there exists a compact contractible
$n$-polyhedron $Q$ such that $Q$ does not embed in $\R^{2n-1}$ but $Q\x Q$ embeds in $\R^{4n-1}$.
\end{maincor}

Theorem \ref{t1} and Corollary \ref{t2} provide a background for the following positive results.

\begin{mainthm} \label{t3}
Let $X$ be a compact $n$-polyhedron and let $m\ge\frac{3(n+1)}2$.
The following are equivalent for each $k>0$:

(i) $X$ embeds in $\R^m$;

(ii) the $k$-fold suspension $\Sigma^k X$ embeds in $\R^{m+k}$;

(iii) the $k$-fold cone $C^k X$ embeds in $\R^{m+k}$;

(iv) $X*T_k$ embeds in $\R^{m+2k}$, where $T_k$ is the join of $k$ copies of the $3$-point set;

(v) $X*Z_k$ embeds in $\R^{m+2k}$, where $Z_k$ is the $(k-1)$-skeleton of the $2k$-simplex.
\end{mainthm}

\begin{mainthm}\label{t4}
Let $X$ be an acyclic compact $n$-polyhedron and let $m\ge\frac{3(n+1)}2$.
The following are equivalent for each $k>0$:

(i) $X$ embeds in $\R^m$;

(ii) $X\x I^k$ embeds in $\R^{m+k}$;

(iii) $X\x\text{(triod)}^k$ embeds in $\R^{m+2k}$.
\end{mainthm}

In Theorem \ref{t3}, (ii)$\Rightarrow$(iii) is obvious and (iii)$\Rightarrow$(i) is immediate%
\footnote{See Corollary \ref{basic-a'} below and \cite{RS}*{Lemma 2.19}.}.
See also Theorem \ref{10.3} below for a simple algebraic proof of (ii)$\Rightarrow$(i).

In Theorem \ref{t4}, (ii)$\Rightarrow$(i) is easy (see Theorem \ref{10.6} below) and was first proved in \cite{SS}.

The real content of Theorems \ref{t3} and \ref{t4} are the remaining assertions.
Their proofs originally appeared in the 2006 preprint \cite{MS}.%
\footnote{To be precise, \cite{MS}*{Corollaries 4.4 and 4.7} are the same as 
Theorems \ref{t3} and \ref{t4} except that they allow $X$ to be an $n$-dimensional compactum 
(rather than a polyhedron) at the cost of imposing a slightly stronger restriction 
$m>\frac{3(n+1)}2$.
But their proofs also yield Theorems \ref{t3} and \ref{t4} if one applies the
classical Haefliger--Weber embeddability criterion for polyhedra (see Theorem \ref{basic-a} below) 
instead of the embeddability criterion for compacta obtained in \cite{MS}.}
The implication (iv)$\Rightarrow$(i) in Theorem \ref{t3} was recently rediscovered by S. Parsa 
\cite{Pa0} (see also \cite{PS}).
Parsa pointed out to the author that the proofs of Theorems \ref{t3} and \ref{t4} in \cite{MS}
are not easy to read as they omit many details.
The present note contains a more detailed exposition of these proofs (at least 4 pages instead of 1 page),
with a few minor errors corrected.

\begin{remark}
Let us note that the dimensional restrictions in Theorem \ref{t4}, (ii)$\Rightarrow$(i) 
cannot be dropped.
Indeed, if $X$ is a non-simply-connected homology $n$-ball (i.e.\ a homology sphere minus an open ball), 
then $X$ does not embed in $\R^n$ (by Seifert--van Kampen), but $X\x I$ embeds in $\R^{n+1}$ since 
every homology sphere bounds a contractible topological manifold (by Kervaire and Freedman; 
see \cite{AG}*{Theorem 0}), whose double has to be the sphere by Seifert--van Kampen and 
by the generalized Poincar\'e conjecture.
\end{remark}

\section{Background and two easy theorems}\label{background}

The following notation will be used throughout the paper.
Given a space $X$, we write $\tilde X=X\x X\but\Delta_X$, where $\Delta_X=\{(x,x)\mid x\in X\}$.
We also write $\bar X=\tilde X/t$, where $t$ is the factor-exchanging involution, $t(x,y)=(y,x)$.
The $m$-sphere $S^m$ will be understood to be endowed with the antipodal involution.
For a $\Z/2$-space $X$ we denote by $\Sigma^k X$ its unreduced $k$-fold suspension $S^{k-1}*X$ with 
the diagonal action of $\Z/2$.

If $X$ embeds in $\R^m$, then there exists an equivariant map $\tilde X\to S^{m-1}$.
Namely, if $g\:X\to\R^m$ is an embedding, then $\tilde g\:\tilde X\to S^{m-1}$ 
is defined by $\tilde g(x,y)=\frac{g(x)-g(y)}{||g(x)-g(y)||}$.

\begin{theorem} \label{basic-a}
{\rm (Haefliger--Weber \cite{We}; see also \cite{Ad}, \cite{M1}*{3.1})}
If $X$ is a compact $n$-polyhedron such that there exists an equivariant map
$\tilde X\to S^{m-1}$ and $m\ge\frac{3(n+1)}2$, then $X$ piecewise linearly embeds in $\R^m$.
\end{theorem}

\begin{corollary} \label{basic-a'} If $X$ is a compact $n$-polyhedron and $m\ge\frac{3(n+1)}2$, 
then $X$ embeds in $\R^m$ if and only if it PL embeds in $\R^m$.
\end{corollary}

\begin{lemma} \label{basic-b}
{\rm (Conner--Floyd \cite{CF}; see also \cite{AM}*{2.7(a)})}
If $K$ is a $k$-polyhedron with a free involution, the suspension map
$[K,S^{m-1}]_{\Z/2}\to[\Sigma K,S^m]_{\Z/2}$ is surjective for $k\le2m-3$ and 
injective for $k\le 2m-4$.
\end{lemma}

\begin{theorem}\label{10.3} Let $X$ be a compact $n$-polyhedron,
$Y$ be a compact polyhedron such that there exists an equivariant map $S^{k-1}\to\tilde Y$,
and assume that $m\ge\frac{3(n+1)}2$.
If $X*Y$ admits a level-preserving (i.e.\ commuting with the projections $X*Y\xr{p} I$ and 
$\R^{m+k}\x I\to I$) embedding in $\R^{m+k}\x I$, then $X$ embeds in $\R^m$.
\end{theorem}

To prove (ii)$\Rightarrow$(i) in Theorem \ref{t3} (in the case $k=1$, which implies the general case) 
we apply the case $Y=S^0$, $k=1$ of Theorem \ref{10.3} to the embedding $e\x p\:X*S^0\to\R^{m+1}\x I$, 
where $e\:X*S^0=\Sigma X\to\R^{m+1}$ is the given embedding.

\begin{proof} The subset $Z$ of $\widetilde{X*Y}$ consisting of all pairs $(a,b)$
such that $p(a)=p(b)$, contains a copy of $\tilde X*\tilde Y$.
So the given embedding yields an equivariant map
$\Sigma^k\tilde X\to\tilde X*\tilde Y\subset Z\to S^{m+k-1}$.
Since $2n\le 2m-3$, by Lemma \ref{basic-b} it desuspends to an equivariant map 
$\tilde X\to S^{m-1}$.
Hence by Theorem \ref{basic-a} $X$ embeds in $\R^m$.
\end{proof}

\begin{theorem}\label{10.6} Let $X$ be an acyclic compact $n$-polyhedron,
$m\ge\frac{3(n+1)}2$.
If $X$ ``can be instantaneously taken off itself in $\R^{m+1}$'', that is,
the mapping cylinder of the projection $X\sqcup X\xr{\pi} X$ admits
a level-preserving (i.e.\ commuting with the projections $MC(\pi)\xr{p} I$ and 
$\R^{m+1}\x I\to I$) embedding in $\R^{m+1}\x I$, then $X$ embeds in $\R^m$.
\end{theorem}

To prove (ii)$\Rightarrow$(i) in Theorem \ref{t4} (in the case $k=1$, which implies the general case)
we apply Theorem \ref{10.6} to the embedding $e\x p\:MC(\pi)\to\R^{m+1}\x I$, where 
$e\:MC(\pi)\cong X\x I\to\R^{m+1}$ is the given embedding.

The idea of the following proof comes from \cite{SS}.

\begin{proof}
Let $g\:X\emb\R^{m+1}$ be the given embedding that can be instantaneously
taken off itself.
Then $\tilde g\:\tilde X\to S^m$ is homotopic to a map $\phi$ that extends to $X\x X$.
Since $X$ is acyclic, $\phi$ is null-homotopic.
Combining the null-homotopy with its reverse, we get an equivariant map
$\Sigma\tilde X\to S^m$.
Since $2n\le 2m-3$, by Lemma \ref{basic-b} it desuspends to an equivariant map 
$\tilde X\to S^{m-1}$.
Hence by Theorem \ref{basic-a} $X$ embeds in $\R^m$.
\end{proof}

\section{Geometric cohomology}

To prove Theorem \ref{t1} we will use the geometric description of cohomology as developed by Buoncristiano, 
Rourke and Sanderson \cite{BRS} (see also \cite{Ro}, \cite{Fe}*{Chapter 1}).
A quick introduction is in \cite{M1}*{\S2}, and we will follow it for terminology.
Here we only hint at the definitions, and then discuss two examples which will be relevant to the proof 
of Theorem \ref{t1}.

Let us first recall that ordinary homology $H_n(P)$ of the polyhedron $P$ with integer coefficients is isomorphic
to the bordism group of {\it singular oriented $n$-pseudomanifolds} in $P$, that is, of maps of the form $f\:M\to P$, 
where $M$ is a (possibly empty) closed oriented $n$-dimensional pseudomanifold (see \cite{M1}*{\S2} for 
a more detailed statement and \cite{Fe}*{Theorem 1.3.7} for the proof).
What we need is the dual description of cohomology.

Let $K$ be a simplicial complex.
An {\it $n$-comanifold} in $K$ (called ``$n$-mock bundle'' in \cite{BRS}) is a PL map of the form $f\:X\to |K|$, 
where $X$ is a polyhedron (possibly empty), 
such that for each $i$-simplex $\sigma$ of $K$, $f^{-1}(\sigma)$ is a compact PL $(i-n)$-manifold with boundary, 
and moreover its boundary coincides with $f^{-1}(\partial\sigma)$.
(If $k<0$, then all $k$-manifolds are empty.)
Thus $n$ is the {\it co}dimension of the comanifold, and we do not care about its dimension.
Let us note that by definition $f$ is transverse to the triangulation $K$ (see \cite{BRS} concerning 
PL transversality).
For instance, if $T$ is the triod $\{0\}*\{1,2,3\}$, then every
embedded $1$-comanifold $G$ in the three-page book $T\x\R$ is a graph with vertices of degrees
$1$ and $3$, where the degree $3$ vertices coincide with the intersection
points of $G$ with the binding $\{0\}\x\R$, and the degree $1$ vertices
coincide with the intersection points of $G$ with the page edges $\{1,2,3\}\x\R$. 

A {\it cobordism} between two $n$-comanifolds $f_0\:X_0\to |K|$ and $f_1\:X_1\to |K|$ is an $n$-comanifold
$\Phi\:W\to|K|\x I$ such that $\Phi^{-1}(|K|\x\{j\})=X_j$ and $\Phi|_{X_j}=f_j$ for $j=0,1$.
It is straightforward to define {\it co-oriented $n$-copseudomanifolds} and their cobordisms.
The cobordism group of co-oriented $n$-copseudomanifolds is isomorphic to the ordinary cohomology $H^n(|K|)$ 
with integer coefficients (see \cite{M1}*{\S2} for a more detailed statement and \cite{BRS} for the proof).
It should also be noted that if $L$ is a subdivision of $K$, then every $n$-comanifold in $L$ is 
an $n$-comanifold in $K$, and every $n$-comanifold in $K$ is cobordant to an $n$-comanifold in $L$; 
the same can be said of co-oriented $n$-copseudomanifolds \cite{BRS}.

\begin{remark}
(a) When $n=1$, the geometric description of $H^1(|K|)$ is closely related to the Pontryagin 
construction for maps $|K|\to S^1$ and their homotopies. 

(b) When $|K|$ is an oriented PL $m$-manifold, co-oriented $n$-copseudomanifolds in $K$ are the same objects as 
singular oriented $(m-n)$-pseudomanifolds in $|K|$ that are transverse to the triangulation $K$ \cite{BRS}, 
and thus the geometric description of $H^n(|K|)$ coincides with the geometric description of 
the Poincar\'e dual group $H_{m-n}(|K|)$ modulo the PL transversality theorem of \cite{BRS}.

(c) If $c$ is a simplicial $n$-cocycle in $K$, the corresponding co-oriented $n$-copseudomanifold $f\:X\to |K|$
can be constructed as follows.
If $\sigma$ is an oriented $n$-simplex of $K$, let $X_\sigma$ be the disjoint union of $|c(\sigma)|$ copies 
$|\sigma^*|_1,\dots,|\sigma^*|_{|c(\sigma)|}$ of the dual cone%
\footnote{$\sigma^*$ is the subcomplex of the barycentric subdivision $K'$ consisting of all simplexes 
of the form $\hat\tau_1*\dots*\hat\tau_n$, where $\sigma\subset\tau_1\subsetneqq\dots\subsetneqq\tau_n$ 
are simplexes of $K$ and $\hat\tau$ denotes the barycenter of $\tau$.
If $K$ is a combinatorial $n$-manifold and $\sigma$ is a $k$-simplex, then $|\sigma^*|$ is an $(n-k)$-cell 
intersecting $\sigma$ at $\hat\sigma$.}
$|\sigma^*|$, co-oriented in $|K|$ according to the orientation of $\sigma$ and the sign of $c(\sigma)$.
If $\tau$ is an oriented $(n+1)$-simplex of $K$, the condition $\delta c(\tau)=0$ implies that if
$\partial\tau=\sigma_1+\dots+\sigma_k$, then $|c(\sigma_1)|+\dots+|c(\sigma_k)|$ is even, and moreover
the co-oriented dual cones in $X_{\sigma_1}\sqcup\dots\sqcup X_{\sigma_k}$ can be matched in pairs, so that
each $|\sigma_i^*|_{i'}$ is matched with some $|\sigma_j^*|_{j'}$, where $j\ne i$, respecting their co-orientations.
Let $X$ be the quotient of $\bigsqcup_\sigma X_\sigma$, where $\sigma$ runs over all $n$-simplexes of $K$,
each with some fixed orientation, by gluing for each oriented $(n+1)$-simplex $\tau$ of $K$ each matched pair
$|\sigma_i^*|_{i'}$, $|\sigma_j^*|_{j'}$ along the dual cone $|\tau^*|$.
Let $f\:X\to |K|$ be the obvious projection, sending the image of each $|\sigma^*|_i$ onto $\sigma^*$.
\end{remark}

\begin{example} \label{e1}
Let $M$ be a non-orientable connected closed PL $n$-manifold.
Let us fix a triangulation of $M$ and an orientation of each $n$-simplex of the triangulation.
Then a co-oriented $n$-copseudomanifold in $M$ is a finite linear combination%
\footnote{If $p$ is a co-oriented point in $M$, then $3p$ stands for the similarly co-oriented map of 
the $3$-point set into $M$, sending all $3$ points onto $p$, and $-p$ stands for $p$ with 
the opposite co-orientation.}
of co-oriented points in the interiors of $n$-simplexes of $M$, where the co-orientations of points 
are induced by the orientations of the $n$-simplexes.
If $p$ is such a co-oriented point in $M$, and $\phi\:I\to M$ is a PL path from $p$ to itself that is orientation-reversing 
(in the sense of local orientation of $M$) and transverse
to the triangulation of $M$, then $\Phi\:I\to M\x I$, $\Phi(t)=\big(\phi(t),t\big)$, is a cobordism from $p$ to $-p$.
If $q$ is another such co-oriented point in $M$, then by connecting $q$ by a path to $p$ we similarly obtain 
a cobordism from $q$ to either $p$ or $-p$.
And in fact every connected cobordism between nonempty $n$-copseudomanifolds is of this form (since $M$ is compact).
Thus we have computed geometrically that $H^n(M)\simeq\Z/2$.
\end{example}

\begin{example} \label{e2} Let $M$ be an oriented connected closed PL $n$-manifold.
Let $N$ be the polyhedron obtained from $M$ by first removing the interior $\mathring D$ of a closed 
PL $n$-ball $D$ and then gluing $D$ back in along some PL map $f\:S^{n-1}\to S^{n-1}$ of degree $2$
which is a ramified covering (for instance, the $(n-2)$-fold suspension over the double covering $S^1\to S^1$).
Let $C$ and $\mathring C$ denote the images of $D$ and $\mathring D$ in $N$ and let $N_0=N\but\mathring C$.
Let us note that both $\mathring C$ and $N_0$ inherit their orientations from that of $M$; moreover, 
at each regular value of the ramified covering $f$, each of the two sheets of $\mathring C$ is attached to $N_0$ 
respecting the orientations (since $\deg f>0$).
Let us fix some triangulation of $N$.
Then a co-oriented $n$-copseudomanifold in $N$ is a finite linear combination of co-oriented points 
in the interiors of simplexes of $N$, where the co-orientations of points are induced by the orientations 
of $\mathring C$ and $N_0$.
If $p$ is such a co-oriented point in the interior of $N_0$, and $r$ is such a co-oriented point in 
$\mathring C$, let us show that $p$ is cobodrant to $2r$.
Indeed, let $T$ be the triod and let $\pi\:T\to I$ be the map sending
the cone apex to $\frac12$, one leaf to $0$ and the other two leaves to $1$, and linear on the three edges.
Let $\phi\:T\to N$ be a PL map sending the cone apex to a point in $\partial N_0$, one leg of the triod, namely, 
$\pi^{-1}([0,\frac12])$, entirely into $N_0$, so that it ends at $p$, and the other two legs entirely into $C$, 
so that they both end at $r$.
If $\phi$ is transverse to the triangulation of $N$, then $\Phi\:T\to N\x I$, $\Phi(t)=\big(\phi(t),\pi(t)\big)$, is
easily seen to be a cobordism between $p$ and $2r$.
Using this, by arguing as in Example \ref{e1} we will compute geometrically that $H^n(N)\simeq\Z$.
\end{example}

\section{Proof of Parsa's theorem} \label{proof}

Let $K$ be an $n$-dimensional join of the $k_i$-skeleta of $(2k_i+2)$-simplexes.
(For instance, the $n$-skeleton of the $(2n+2)$-simplex or the join of $n+1$ copies of the $3$-point set.)

\begin{remark}
To motivate what follows let us recall that $K$ is the simplest example of a polyhedron with non-zero 
van Kampen obstruction (see \cite{M1}*{Examples 3.3, 3.5}).
Moreover, by removing the interior of an $n$-disk $D$ lying in the interior of some $n$-simplex of $K$
and then gluing $D$ back in along some PL map $S^{n-1}\to S^{n-1}$ of degree $2$ we obtain a polyhedron
which also has non-zero van Kampen obstruction, but zero mod~$2$ van Kampen obstruction \cite{M1}*{Example 3.6}.
\end{remark}

Let us fix two disjoint $n$-simplexes $A$ and $B$ in the natural triangulation of $K$.
Let $p$ be a point in the interior of $A$ and $q$ be a point in the interior of $B$.
Let $\pi\:\tilde K\to\bar K$ be the quotient map (see notation in \S\ref{background}).
Let us fix some triangulation of $\bar K$ (it has to be infinite since $\bar K$ is non-compact).
The unordered pair $\{p,q\}$ may be identified with a point of $\bar K$.
We may assume that it lies in the interior of some $2n$-simplex of $\bar K$.
It is not hard to construct an explicit generic PL map $f\:K\to\R^{2n}$ with double point 
$f(p)=f(q)$ and no other double points (see \cite{M1}*{Example 3.5}).
Therefore the van Kampen obstruction $\theta(K)\in H^{2n}(\bar K)$ 
(see \cite{M1}*{\S3, subsection ``Geometric definition of $\theta(X)$''}) is represented 
by the point $\{p,q\}$ with some co-orientation (which is regarded as a co-oriented 
$2n$-copseudomanifold in $\bar K$).

Let us fix this co-orientation of $\{p,q\}$.
Let us also fix some co-orientations of $p$ and $q$ so that $\{p,q\}=\pi(p,q)$ 
as co-oriented points, where the co-orientation of $(p,q)$ is understood to be 
induced by those of $p$ and $q$, and the co-orientation of $\pi(p,q)$ is understood to be 
induced by that of $(p,q)$ via the double covering $\pi$.
Since $\theta(K)$ is the image of the generator of $H^{2n}(\R P^\infty)\simeq\Z/2$ under the
homomorphism induced by the map $\bar K\to\R P^\infty$ classifying the double cover
$\tilde K\to\bar K$ (see \cite{M1}*{\S3, subsection ``Van Kampen obstruction''}), we have
$\theta(K)=-\theta(K)$.
Hence there exists a cobordism $W_0$ in $\bar K$ between $\{p,q\}$ and $-\{p,q\}$.
It lifts to a cobordism $W$ in $\tilde K$ between $(p,q)$ and $-t(p,q)$, using that 
$\pi\big(t(p,q)\big)=\pi(p,q)=\{p,q\}$ as co-oriented points.
(Here the co-orientation of $t(p,q)$ is understood to be induced from the co-orientation of $(p,q)$ 
via the homemorphism $t$ and has a priori no relation with the co-orientation of $(q,p)$, which is 
induced directly from the co-orientations of $p$ and $q$.
In fact, it is not hard to see that $t(p,q)=(-1)^n(q,p)$ as co-oriented points, but we do not need this.)%
\footnote{One of the referees prefers an alternative explanation. 
According to him/her, the oriented point of $\bar K$ that was denoted $\{p,q\}$ is ``in reality''
$\{(p,q), (-1)^n (q,p)\}$. 
Then $-\{p,q\} = \{-(p,q), (-1)^{n+1} (q,p)\}$, and the cobordism $W$ matches the orientation of $(p,q)$ 
with that of $(-1)^{n+1}(q,p)=-t(p,q)$.}

Let us pick a point $x$ in the interior of $A$ such that $\{x\}\x K\x I$ is disjoint from 
the image of $W$.
Then $x$ lies in a closed $n$-ball $D$ contained in the interior of $A$ such that $D\x K\x I$ is 
still disjoint from the image of $W$.
In particular, we have $p\in A\but D$.
Let $L$ be obtained from $K$ by first removing the interior $\mathring D$ of $D$ and then gluing $D$ 
back in along some PL map $f\:S^{n-1}\to S^{n-1}$ of degree $2$ which is a ramified covering.
(Thus $f$ must be the $2$-fold covering when $n=2$.)
Let $C$ and $\mathring C$ be the images of $D$ and $\mathring D$ in $L$ (thus $C\cong\R P^2$ when $n=2$), 
and let $A'\subset L$ be the effect of the modification on $A\subset K$.
Then there still exists the cobordism $W$ in $\tilde L$ between $(p,q)$ and $-t(p,q)$, and 
by the following lemma $\theta(L)\in H^{2n}(\bar L)$ is still represented by the co-oriented point 
$\{p,q\}$ of $\bar L$.

\begin{lemma} The homeomorphism $h\:C\but\mathring C\to\partial D=\partial D\x pt$ extends to an embedding 
$C\to D\x I^n$.
\end{lemma}

\begin{proof} Let us extend $h$ to a generic PL map $f\:C\to D\x I^n$.
Then $f$ has only isolated double points in $\mathring C$.
Let $f(x)=f(y)$ be one.
Let us connect $x$ and $y$ by an arc $J$ in $\mathring C$, disjoint from the preimages of
the other double points (using that $n>1$).
If $n>2$, then $f(J)$ bounds a $2$-disk $\Delta$ in $D\x I^n$ which meets $f(C)$ only in 
$\partial\Delta$.
Then a small regular neighborhood $R$ of $\Delta$ (which is a $2n$-ball) meets $f(C)$ 
in the image of a small regular neighborhood $N$ of $J$ (which is an $n$-ball).
Now $f$ embeds $\partial N$ in $\partial R$, and we redefine $f|_N\:N\to R$ to be
the conical extension of that embedding.
By repeating the same construction for each double point of $f$ we obtain the desired 
embedding for $n>2$.
When $n=2$, we obtain it by using that $\R P^2$ embeds in $\R^4$ and $S^1$ unknots in $\R^4$.
\end{proof}

Let $r$ be a point in $\mathring C$.
Similarly to Example \ref{e2} there exists a cobordism $V_0$ in $A'$ between $p$ and $2r$, 
for an appropriate co-orientation of $r$.
Since $A'\x B$ is a neighborhood of $V_0\x\{q\}$ in $\tilde L$, $V:=V_0\x\{q\}$ is a cobordism
in $\tilde L$ between $(p,q)$ and $2(r,q)$.

Now let us consider the co-oriented point $\big((p,q),(p,q)\big)$ in $\tilde L\x\tilde L$.
Since is $\tilde L\x A'\x B$ is a neighborhood of $V\x\{(p,q)\}$ in $\tilde L\x\tilde L$,
$V\x\{(p,q)\}$ is a cobordism in $\tilde L\x\tilde L$ between $\big((p,q),(p,q)\big)$ and
$2\big((r,q),(p,q)\big)$.
In a similar way we obtain the following sequence of cobordisms in $\tilde L\x\tilde L$:
\begin{multline*}
\scalebox{0.9}{$\big((r,q),(p,q)\big)\xR{\{(r,q)\}\x W}
-\big((r,q),t(p,q)\big)\xR{\{(r,q)\}\x tV}
-2\big((r,q),t(r,q)\big)\xL{V\x\{t(r,q)\}}
-\big((p,q),t(r,q)\big)$}\\
\scalebox{0.9}{$\xR{W\x\{t(r,q)\}}
\big(t(p,q),t(r,q)\big)\xR{tV\x\{t(r,q)\}}
2\big(t(r,q),t(r,q)\big)\xL{\{t(r,q)\}\x tV}
\big(t(r,q),t(p,q)\big)$.}
\end{multline*}
Since $\tilde L\x\tilde L$ can be identified 
via $\big((a,b),(c,d)\big)\mapsto\big((a,c),(b,d)\big)$ with an open subset of 
$\widetilde{L\x L}$, all these cobordisms also take 
place in $\widetilde{L\x L}$ with an appropriate triangulation.

From the usual inclusion $L\x L\subset L*L$ we have $\widetilde{L\x L}\subset\widetilde{L*L}$.
Let us fix some co-orientation of $\widetilde{L\x L}$ in $\widetilde{L*L}$ and also 
a triangulation of $\widetilde{L*L}$ such that $\widetilde{L\x L}$ is transversal to it.
Then every co-oriented point $\big((a,b),(c,d)\big)$ in 
$\tilde L\x\tilde L\subset\widetilde{L\x L}$ which lies in the interior of a $(4n+2)$-simplex 
of $\widetilde{L*L}$ becomes co-oriented in $\widetilde{L*L}$,
and the factor exchanging involution $T$ of $\widetilde{L*L}$ acts on such points by
$T\big((a,b),(c,d)\big)=-\big(t(a,b),t(c,d)\big)$.
Consequently we have the following cobordisms in $\widetilde{L*L}$:
\[\big((p,q),(p,q)\big)\Longrightarrow 2\big((r,q),(p,q)\big)\text{\quad and\quad}
\big((r,q),(p,q)\big)\Longrightarrow -T\big((r,q),(p,q)\big).\]
Using the double cover $\widetilde{L*L}\to\overline{L*L}$ to induce co-orientations
on points of $\overline{L*L}$, we get the following cobordisms in $\overline{L*L}$:
\[\big[(p,q),(p,q)\big]\Longrightarrow 2\big[(r,q),(p,q)\big]\text{\quad and\quad}
\big[(r,q),(p,q)\big]\Longrightarrow -\big[(r,q),(p,q)\big],\]
where $\big[(a,b),(c,d)\big]=\big\{(a,c),(b,d)\big\}\in\overline{L*L}$.
Thus the cobordism class of $\big[(p,q),(p,q)\big]$ is trivial.
But by Proposition \ref{map-join} $\theta(L*L)\in H^{4n+2}(\overline{L*L})$
is represented by $\pm\big[(p,q),(p,q)\big]$.
Thus $\theta(L*L)=0$.
Hence by the Shapiro--Wu theorem (see \cite{M1}*{Theorem 3.2}) $L*L$ embeds 
in $\R^{4n+2}$.

\begin{proposition} \label{map-join}
Let $f\:P^m\to\R^{2m}$ and $g\:Q^n\to\R^{2n}$ be generic PL maps of 
compact polyhedra, and let $\Delta_f=\{(x,y)\in\tilde P\mid f(x)=f(y)\}$.
For any two equivariant maps $\phi\:\Delta_f\to S^0$ and $\psi\:\Delta_g\to S^0$
there exists a generic PL map $h\:P*Q\to\R^{2(m+n+1)}$ such that 
$\Delta_h$ is the subset of $\Delta_f\x\Delta_g$ consisting of all 
$\big((a,b),(c,d)\big)$ such that $\phi(a,b)=\psi(c,d)$.
\end{proposition}

Here $\Delta_f\x\Delta_g$ is regarded as a subset of 
$P^2\x Q^2\subset P^2*Q^2\subset(P*Q)^2$.

\begin{proof} By identifying $\R^{2m}$ and $\R^{2n}$ with a pair of skew affine 
subspaces in $\R^{2m+2n+1}$, say $\R^{2m}\x\{0\}\x\{1\}$ and $\{0\}\x\R^{2n}\x\{-1\}$, 
we get the map $f*g\:P*Q\to\R^{2m+2n+1}$.
Clearly, $\Delta_{f*g}=
\Delta_f*\Delta_g\cup(\Delta_f*\Delta_Q\but\Delta_Q)\cup(\Delta_P*\Delta_g\but\Delta_P)$ 
as subsets of $P^2*Q^2\subset(P*Q)^2$.
On the other hand, let $\bar\phi\:P\to\R$ and $\bar\psi\:Q\to\R$ be some PL extensions 
of the compositions
$\pi_1(\Delta_f)\xr{(\pi_1|_{\Delta_f})^{-1}}\Delta_f\xr{\phi}S^0$ and 
$\pi_2(\Delta_g)\xr{(\pi_2|_{\Delta_g})^{-1}}\Delta_g\xr{\psi}S^0$, where $S^0=\{-1,1\}$
and $\pi_i\:P\x P\to P$ and $\pi_i\:Q\x Q\to Q$ denote the projections onto 
the $i$th factor.
Let $\Xi\:P*Q\to\R$ be the sum of the compositions
$\Phi\:P*Q\xr{\bar\phi*\const}\R*\{r\}\xr{\chi}\R$ and 
$\Psi\:P*Q\xr{\const*\bar\psi}\{r\}*\R\xr{\chi}\R$, where $\chi\big(tx+(1-t)r\big)=tx$.
That is, $\Xi(x)=\Phi(x)+\Psi(x)$.
Finally, let $h$ be the joint map
$(f*g)\x\Xi\:P*Q\to\R^{2m+2n+1}\x\R$.
Thus $\Delta_h=\Delta_{f*g}\cap\Delta_{\Xi}$.

Each $x\in\Delta_f*\Delta_Q\but\Delta_Q$ is of the form $x=t(a,b)+(1-t)(c,c)$ for some 
$(a,b)\in\Delta_f$, $c\in Q$ and $t>0$.
We may also write $x=(y,z)$, where $y=ta+(1-t)c$ and $z=tb+(1-t)c$.
We have $\Phi(y)=t\bar\phi(a)=t\phi\big((\pi_1|_{\Delta_f})^{-1}(a)\big)=t\phi(a,b)$ and similarly $\Phi(z)=t\phi(b,a)$.
Since $\phi$ is equivariant and $t>0$, we have $\Phi(y)\ne\Phi(z)$.
On the other hand, $\Psi(y)=\Psi(z)=(1-t)\bar\psi(c)$.
Hence $\Xi(y)\ne\Xi(z)$.
Thus $x\notin\Delta_\Xi$.
Similarly, no point of $\Delta_P*\Delta_g\but\Delta_P$ lies in $\Delta_\Xi$.

Each $x\in\Delta_f*\Delta_g$ is of the form $x=t(a,b)+(1-t)(c,d)$ for some $(a,b)\in\Delta_f$
and $(c,d)\in\Delta_g$.
We may also write $x=(y,z)$, where $y=ta+(1-t)c$ and $z=tb+(1-t)d$.
We have $\Phi(y)=t\bar\phi(a)=t\phi(a,b)$ and similarly $\Phi(z)=t\phi(b,a)$.
Likewise $\Psi(y)=(1-t)\psi(d,c)$ and $\Psi(z)=(1-t)\psi(c,d)$.
Since $\phi$ and $\psi$ are equivariant, we have $\Phi(y)-\Phi(z)=2t\phi(a,b)$
and $\Psi(y)-\Psi(z)=2(1-t)\psi(d,c)=-2(1-t)\psi(c,d)$.
Hence $\Xi(y)-\Xi(z)=2t\phi(a,b)-2(1-t)\psi(c,d)$.
Thus $x\in\Delta_\Xi$ if and only if $t=\frac12$ and $\phi(a,b)=\psi(c,d)$.
Moreover, it is clear that every such double point is transverse.
\end{proof}

\section{Extraordinary van Kampen obstruction}

We will give two slightly different proofs of Theorems \ref{t3} and \ref{t4}.
The first approach does {\it not} use the definitions of the present section, and 
the second approach has some other advantages.%
\footnote{In particular, the proof of Proposition \ref{10.4}(a) based on the extraordinary
van Kampen obstruction is considerably simpler than its alternative proof which does not use it.
Also, this alternative proof would not easily generalize to the case of compacta.}

For a space $Q$ let $Q_+$ denote the pointed space $Q\sqcup pt$ with basepoint at $pt$.
If $Q$ is a free $\Z/2$-space, for instance, $Q=\tilde K$ for some polyhedron $K$,
then $Q_+$ is a pointed $\Z/2$-space.

If $G$ is a finite group and $V$ is a finite-dimensional $\R G$-module, let $S^V$ denote 
the one-point compactification of the Euclidean space $V$ with the obvious action of $G$.
Let $mT$ denote the sum of $m$ copies of the nontrivial one-dimensional real representation $T$ of $\Z/2$.
Thus $S^{mT}$ is the $m$-sphere with the action of $\Z/2$ fixing the basepoint at infinity 
and restricting to the sign action $x\mapsto -x$ on $\R^m$.

Let $P$ be a pointed polyhedron with an action of $G$ that fixes the basepoint $*$, and let us assume that
$P$ is $G$-homotopy equivalent to a compact polyhedron.
Then the equivariant stable cohomotopy group
$$\omega^{V-W}_G(P):=[S^{W+V_\infty}\wedge P,\ S^{V+V_\infty}]^*_G$$ 
is well-defined, where $V_\infty$ denotes a sufficiently large (with respect to the partial ordering with respect 
to inclusion) finite-dimensional $\R G$-submodule of the countable direct sum $\R G\oplus\R G\oplus\dots$
(see \cite{M+}, \cite{Ada}).
Let us note that the stabilization map $[\tilde K_+,S^{mT}]^*_{\Z/2}\to\omega^{mT}_{\Z/2}(\tilde K_+)$
is bijective for $m\ge\dim K+1$ (Bredon--Hauschild; see \cite{M+}*{IX.I.4}, \cite{Ada}*{3.3}).

The following {\it extraordinary van Kampen obstruction}
$$\Theta^m(\tilde K)\in\omega^{mT}_{\Z/2}(\tilde K_+)$$
to embeddability of a compact polyhedron $K$ in $\R^m$ was introduced in 
\cite{M1}*{\S6} (where it is denoted $\Theta^m(K)$).
Let $Q$ be a polyhedron with a free action of $\Z/2$ such that $Q$ is equivariantly homotopy equivalent 
to a compact polyhedron.
Let $f\:Q\to S^\infty$ be an arbitrary equivariant map (with respect to the antipodal involution on $S^\infty$),
and let $\alpha\in [Q_+,S^{mT}]^*_{\Z/2}$ be the class of the composition
$$Q\xr{f} S^\infty\to S^\infty/S^{m-1}\cong S^\infty_+\wedge S^{mT}\to S^{mT}.$$
Then $\Theta^m(Q)\in\omega^{mT}_{\Z/2}(Q_+)$ is defined to be the image of $\alpha$ under the stabilization map.

If $K$ is a compact polyhedron that embeds in $\R^m$, then any embedding $g\:K\to\R^m$ yields an equivariant map 
$f\:\tilde K\xr{\tilde g} S^{m-1}\subset S^\infty$, and hence $\Theta^m(\tilde K)=0$.

\section{Proofs of stability results}

Let $X$ and $Y$ be compacta.
Given a $y\in Y$ and an $s\in I=[0,1]$, let $[y,s]$ denote the image of
$(y,s)$ in $CY=Y\x I\,\big/\,Y\x\{0\}$.
Given additionally an $x\in X$, let $[x,y,s]$ denote the image of $(x,y,s)$ in 
$X*Y=X\x Y\x I\,\big/\,\pi$, where $\pi\:X\x Y\x\partial I\to X\sqcup Y$ projects 
$X\x Y\x\{0\}$ onto the first factor and $X\x Y\x\{1\}$ onto the second factor.
Thus $[x,y,0]$ is independent of $y$, and $[x,y,1]$ is independent of $x$.

Let $Y^\urcorner$ be the subset of $\widetilde{CY}$ consisting of all pairs 
$\big([y,s],[y',s']\big)$ such that either $s=1$ or $s'=1$, and if $s<1$ or $s'<1$, 
then $y\ne y'$.
It is easy to see that $Y^\urcorner$ is equivariantly homeomorphic to the double 
mapping cylinder of the projections $Y\xleftarrow{p_1}\tilde Y\xr{p_2}Y$.

\begin{lemma}\label{cone-dp} $\widetilde{CY}$ is equivariantly homotopy equivalent 
to $Y^\urcorner$.
\end{lemma}

\begin{proof}
First we observe that $\widetilde{CY}$ equivariantly deformation retracts onto its
intersection $Z$ with $Y\x CY\cup CY\x Y$.
Indeed, $Z$ consists of all pairs $\big([y,s],[y',s']\big)$ such that either 
$s=1$ or $s'=1$, and if $s=s'=1$, then $y\ne y'$.
A homotopy $r_t\:\widetilde{CY}\to\widetilde{CY}$ is given by $r_t\big([y,s],[y',s']\big)=
\big([y,(1-t)s+t\min(\frac{s}{s'},1)],\,[y',(1-t)s'+t\min(\frac{s'}{s},1)]\big)$.
Here $\min(\frac{\sigma}0,1)$ is understood to be $1$ if $\sigma>0$ and we never have 
$s=s'=0$.
The homotopy is well-defined since $(1-t)\sigma+t\min(\frac{\sigma}{\sigma'},1)=0$ 
if $\sigma=0$ and $\sigma'>0$.
It is easy to check that $r_t\big([y,s],[y',s']\big)\in\widetilde{CY}$ for each $t\in I$.%
\footnote{Indeed, if $y=y'$, then $s\ne s'$.
Let us show that $(1-t)s+t\min(\frac{s}{s'},1)\ne(1-t)s'+t\min(\frac{s'}{s},1)$.
If $\sigma>\sigma'$, then $(1-t)(\sigma-\sigma')+t\frac{\sigma-\sigma'}{\sigma}=
(\sigma-\sigma')\big(1+t(\frac{1}{\sigma}-1)\big)>0$ 
and consequently $(1-t)\sigma+t>(1-t)\sigma'+t\frac{\sigma'}{\sigma}$.}
It follows that $r_t$ is an equivariant deformation retraction of $\widetilde{CY}$ onto $Z$.

Finally, $Z$ is equivariantly homotopy equivalent to $Y^\urcorner$.
Indeed, let $h_t\:Z\to Z$ be given by 
$h_t\big([y,s],[y',s']\big)=\big([y,\,td(s)+(1-t)s],[y',\,td(s')+(1-t)s']\big)$, 
where $d(\sigma)=\lfloor\sigma\rfloor$ (that is, $1$ if $\sigma=1$ and $0$ if $\sigma<1$) if $y=y'$
and $d(\sigma)=1-\min\big\{\frac{1-\sigma}{\min\{d(y,y'),\,1\}},\,1\big\}$ if $y\ne y'$.
Clearly $h_t$ is continuous and equivariant, $h_1(Z)\subset Y^\urcorner$ and 
$h_t(Y^\urcorner)\subset Y^\urcorner$ for all $t\in I$.
Then $\id_Z$ is homotopic the composition $Z\xr{h_1}Y^\urcorner\subset Z$ to via $h_t$, and $\id_{Y^\urcorner}$ 
is homotopic to the composition $Y^\urcorner\subset Z\xr{h_1}Y^\urcorner$ via $h_t|_{Y^\urcorner}$.
\end{proof}

\begin{proposition}\label{10.1} Let $X$ be a compact $n$-polyhedron and $Y$ be 
a compact polyhedron of dimension $\le\frac {m-n+k-3}2$ such that there exists 
an equivariant map $S^{k-1}\to\widetilde{CY}$.

(a) If $n\le m-2$ and there exists an equivariant map $\widetilde{X*Y}\to S^{m+k-1}$,
then there exists an equivariant map $\tilde X\to S^{m-1}$.

(b) If $\Theta^{m+k}(\widetilde{X*Y})=0$, then $\Theta^m(\tilde X)=0$.
\end{proposition}

Part (a) follows from (b) and \cite{M1}*{Lemma 6.4}.
We also include an alternative proof of (a) which does not use the extraordinary van Kampen 
obstruction.

\begin{proof} 
The beginning of the proof is common for (a) and (b).
By Lemma \ref{cone-dp} $\widetilde{CY}$ is equivariantly homotopy equivalent to $Y^\urcorner$, 
so the hypothesis yields an equivariant map $S^{k-1}\to Y^\urcorner$.
Thus we obtain an equivariant map $\phi\:\Sigma^k\tilde X\to Y^\urcorner*\tilde X$.

Let $Z$ be the subset of $\widetilde{X*Y}$ that is the image of $\tilde X\x\tilde Y\x I^2$
under the quotient map $X^2\x Y^2\x I^2\to (X*Y)^2$; 
in other words, $Z$ consists of all pairs $\big([x,y,s],[x',y',s']\big)$ satisfying
the following two conditions: 
\begin{enumerate}
\item either $x\ne x'$ or $\max(s,s')=1$; and 
\item either $y\ne y'$ or $\min(s,s')=0$.
\end{enumerate}
Let us define an equivariant map $f\:Z\to\tilde X*Y^\urcorner$ by the formula
\[f\big([x,y,s],[x',y',s']\big)=\big[(x,x'),\,\big([y,\min(\tfrac{s}{s'},1)],\,
[y',\min(\tfrac{s'}{s},1)]\big),\,\max(s,s')\big].\]
Here $\min(\frac{\sigma}0,1)$ is understood to be $1$ if $\sigma>0$ and is undefined when $\sigma=0$. 
Let us check that $f$ is well-defined.
Firstly, if $\frac00$ occurs in the right hand side, then $s=0$ and $s'=0$, so $\max(s,s')=0$
and therefore the term involving $\frac00$ need not be defined.
Secondly, if $x$ is undefined in the left hand side, then $s=1$, so $\max(s,s')=1$ and therefore 
$x$ need not be defined in the right hand side.
Similarly for $x'$.
Finally, if $y$ is undefined in the left hand side, then $s=0$ and we consider two cases.
If $s'=0$, then $\max(s,s')=0$ and hence $y$ need not be defined in the right hand side.
If $s'>0$, then $\min(\tfrac{s}{s'},1)=0$ and hence again $y$ need not be defined in the right hand side.
Similarly for $y'$.
It is easy to see that $f$ is continuous and equivariant.

The restrictions of $f$ corresponding to various values of $s$ and $s'$ are as follows:
\[\begin{matrix}
s=s'=0 &\text{a homeomorphism between two copies of }\tilde X\\
\min(s,s')=0,\ \max(s,s')\in(0,1) &\text{a homeomorphism between two copies of }\tilde X\x Y\\
s,s'\in(0,1) &\text{a homeomorphism between two copies of }\tilde X\x\tilde Y\\
\min(s,s')=0,\ \max(s,s')=1 & \text{a copy of the projection }X\x Y\to Y\\
\min(s,s')\in(0,1),\ \max(s,s')=1 &\text{a copy of the projection }X\x\tilde Y\to\tilde Y\\
s=s'=1 &\text{a homeomorphism between two copies of }\tilde Y.
\end{matrix}\]
Thus $f$ is surjective, and its restriction 
$f|_{\dots}\:Z\but f^{-1}(Y^\urcorner)\to\tilde X*Y^\urcorner\but Y^\urcorner$, which 
corresponds to the case $\max(s,s')<1$ (i.e.\ the first three cases of the previous list), 
is a homeomorphism.
Consequently, if $g$ denotes the restriction $f|_{\dots}\:f^{-1}(Y^\urcorner)\to Y^\urcorner$ of $f$,
which corresponds to the case $\max(s,s')=1$ (i.e.\ the last three cases of the previous list),
then the mapping cylinder $MC(f)$ collapses onto $Z\cup MC(g)$.%
\footnote{In fact, $f^{-1}(Y^\urcorner)$ is equivariantly homeomorphic to $X\x Y^\urcorner\,\big/\,p$, 
where $p\:X\x\tilde Y\to\tilde Y$ is the projection, and $g$ is a copy of the map
$X\x Y^\urcorner\,\big/\,p\to Y^\urcorner$ determined by the projection 
$X\x Y^\urcorner\to Y^\urcorner$.}
Since $\dim Y^\urcorner\le m-n+k-2$ and $\dim X=n$, we have $\dim MC(g)\le m+k-1$.

{\it Proof of (b).}
We have $\omega^V_{\Z/2}\big(MC(f),Z\big)\simeq\omega^V_{\Z/2}\big(MC(g)\cup Z,\,Z\big)=0$ 
if $\dim V\ge m+k$.
Hence the restriction map 
$\tilde\omega^{(m+k)T}_{\Z/2}\big(MC(f)_+\big)\to\tilde\omega^{(m+k)T}_{\Z/2}(Z_+)$ is 
an isomorphism.
Therefore so is 
$f^*\:\tilde\omega^{(m+k)T}_{\Z/2}(\tilde X*Y^\urcorner_+)\to\tilde\omega^{(m+k)T}_{\Z/2}(Z_+)$.
Since $\Theta^{m+k}(\widetilde{X*Y})=0$ and $Z\subset\widetilde{X*Y}$, we have 
$\Theta^{m+k}(Z)=0$.
Since $f^*$ is injective, we get that $\Theta^{m+k}(\tilde X*Y^\urcorner)=0$.
The existence of the equivariant map $\phi\:\Sigma^k\tilde X\to Y^\urcorner*\tilde X$
now yields that $\Theta^{m+k}(\Sigma^k\tilde X)=0$.
Hence $\Theta^m(\tilde X)=0$.

{\it Proof of (a).}
Since $\dim MC(g)\le m+k-1$, the given equivariant map $\widetilde{X*Y}\to S^{m+k-1}$ 
restricted to $Z$ extends to an equivariant map $Z\cup MC(g)\to S^{m+k-1}$.
Since $MC(f)$ equivariantly deformation retracts onto $Z\cup MC(g)$, we get
an equivariant map $\psi\:MC(f)\to S^{m+k-1}$.
Then we get the following composition of equivariant maps:
$\Sigma^k\tilde X\xr{\phi}Y^\urcorner*\tilde X\subset MC(f)\xr{\psi}S^{m+k-1}$.
Since $2n\le 2m-3$, by Lemma \ref{basic-b} it desuspends to an equivariant map 
$\tilde X\to S^{m-1}$.
\end{proof}

\begin{proof}[Proof of Theorem \ref{t3}] 
The implications (ii)$\Rightarrow$(i) and (iii)$\Rightarrow$(i) follow from Theorem \ref{basic-a} 
by a {\it repeated} application of Proposition \ref{10.1}(a).
The implications (iv)$\Rightarrow$(i) and (v)$\Rightarrow$(i) follow from 
Proposition \ref{10.1}(a) and Theorem \ref{basic-a} since $\widetilde{CT_k}$ and 
$\widetilde{CZ_k}$ contain $S^{2k-1}$ by the Flores construction (see \cite{M2}). 
\end{proof}

By similar arguments, to prove Theorem \ref{t4} it suffices to prove part (a) of the following

\begin{proposition}\label{10.4} Let $X$ be an acyclic compact $n$-polyhedron and let $Y$ be 
a compact polyhedron such that there exists an equivariant map $S^{k-1}\to\widetilde{CY}$.

(a) If $n\le m-2$ and there exists an equivariant map 
$\widetilde{X\x CY}\to S^{m+k-1}$, then there exists an equivariant map 
$\tilde X\to S^{m-1}$.

(b) If $\Theta^{m+k}(\widetilde{X\x CY})=0$, then $\Theta^m(\tilde X)=0$.
\end{proposition}

Part (a) follows from (b) and \cite{M1}*{Lemma 6.4}.
We also include an alternative proof of (a) which does not use the extraordinary van Kampen 
obstruction.

\begin{proof} The beginning of the proof is common for (a) and (b) and is largely 
similar to that of Proposition \ref{10.1}.
By Lemma \ref{cone-dp} $\widetilde{CY}$ is equivariantly homotopy equivalent 
to $Y^\urcorner$, so the hypothesis yields an equivariant map $S^{k-1}\to Y^\urcorner$.
Thus we obtain an equivariant map $\phi\:\Sigma^k\tilde X\to Y^\urcorner*\tilde X$.

Let $Z$ be the subset of $\widetilde{X\x CY}$ that is the image of
$\tilde X\x\tilde Y\x I\x I\cup X\x X\x\tilde Y\x (I\x\{1\}\cup\{1\}\x I)$
under the quotient map $X^2\x Y^2\x I^2\to(X\x CY)^2$; in other words, $Z$
consists of all pairs $\big((x,[y,s]),(x',[y',s'])\big)$ satisfying
the following two conditions:
\begin{enumerate}
\item either $x\ne x'$ or $\max(s,s')=1$; and 
\item either $y\ne y'$ or $\min(s,s')=0$.
\end{enumerate}
Let us define an equivariant map $f\:Z\to\tilde X*Y^\urcorner$ by 
\[f\big((x,[y,s]),(x',[y',s'])\big)=\big[(x,x'),\,\big([y,\min(\tfrac{s}{s'},1)],\,
[y',\min(\tfrac{s'}{s},1)]\big),\,\max(s,s')\big].\]
Here $\min(\frac{\sigma}0,1)$ is understood to be $1$ if $\sigma>0$ and is undefined when $\sigma=0$. 
Let us check that $f$ is well-defined.
Firstly, if $\frac00$ occurs in the right hand side, then $s=0$ and $s'=0$, so $\max(s,s')=0$
and therefore the term involving $\frac00$ need not be defined.
Finally, if $y$ is undefined in the left hand side, then $s=0$ and we consider two cases.
If $s'=0$, then $\max(s,s')=0$ and hence $y$ need not be defined in the right hand side.
If $s'>0$, then $\min(\tfrac{s}{s'},1)=0$ and hence again $y$ need not be defined in the right hand side.
Similarly for $y'$.
It is easy to see that $f$ is continuous and equivariant.

The restrictions of $f$ corresponding to various values of $s$ and $s'$ are as follows:
\[\begin{matrix}
s=s'=0 &\text{a homeomorphism between two copies of }\tilde X\\
\min(s,s')=0,\ \max(s,s')\in(0,1) &\text{a homeomorphism between two copies of }\tilde X\x Y\\
s,s'\in(0,1) &\text{a homeomorphism between two copies of }\tilde X\x\tilde Y\\
\min(s,s')=0,\ \max(s,s')=1 & \text{a copy of the projection }X^2\x Y\to Y\\
s,s'\in(0,1],\ \max(s,s')=1 &\text{a copy of the projection }X^2\x\tilde Y\to\tilde Y.
\end{matrix}\]

Thus $f$ is surjective, and its restriction 
$f|_{\dots}\:Z\but f^{-1}(Y^\urcorner)\to\tilde X*Y^\urcorner\but Y^\urcorner$, which corresponds
to the case $\max(s,s')<1$ (i.e.\ the first three cases of the previous list), is a homeomorphism.
Consequently, if $g$ denotes the restriction $f|_{\dots}\:f^{-1}(Y^\urcorner)\to Y^\urcorner$ of $f$,
which corresponds to the case $\max(s,s')=1$ (i.e.\ the last two cases of the previous list), then
the mapping cylinder $MC(f)$ equivariantly deformation retracts onto $Z\cup MC(g)$.
In fact, $f^{-1}(Y^\urcorner)$ is equivariantly homeomorphic to $X^2\x Y^\urcorner$ 
and $g$ is a copy of the projection $p\:X^2\x Y^\urcorner\to Y^\urcorner$.
Since $X$ is acyclic, so is $X^2$.

{\it Proof of (b).}
$\omega^V_{\Z/2}\big(MC(f),Z\big)\simeq\omega^V_{\Z/2}\big(MC(g)\cup Z,\,Z\big)
\simeq\omega^V_{\Z/2}\big(MC(g),\,f^{-1}(Y^\urcorner)\big)\simeq
\omega^V_{\Z/2}\big(CX^2\x Y^\urcorner,\,X^2\x Y^\urcorner\big)$ for any $V$.
By map excision the latter group is isomorphic to
$\omega^V_{\Z/2}\big((\{a,b\}*X^2)\x Y^\urcorner,\,\{a\}\x Y^\urcorner\big)$,
where $\Z/2$ acts trivially on $\{a,b\}$.
Since $X^2$ is acyclic, $\{a,b\}*X^2$ is contractible, and hence
$(\{a,b\}*X^2)\x Y^\urcorner$ deformation retracts onto $\{a\}\x Y^\urcorner$.
Thus we get that $\omega^V_{\Z/2}\big(MC(f),Z\big)=0$.
The remainder of the proof is entirely similar to that of Proposition \ref{10.1}(b).

{\it Proof of (a).}
Since $X^2$ is acyclic, $p\:X^2\x Y^\urcorner\to Y^\urcorner$ induces isomorphisms on 
ordinary homology groups by the Vietoris--Begle theorem.
So does $\Sigma^2p\:\Sigma^2(X^2\x Y^\urcorner)\to\Sigma^2 Y^\urcorner$.
Since $\Sigma^2(X^2\x Y^\urcorner)$ and $\Sigma^2 Y^\urcorner$ are simply-connected,
by the Hurewicz theorem $\Sigma^2p$ induces isomorphisms on homotopy groups.
Since $\Sigma^2p$ is equivariant, it descends to a map
$q\:A\to B$, where $A=\Sigma^2(X^2\x Y^\urcorner)/t$ and $B=\Sigma^2 Y^\urcorner/t$.
Clearly $q$ also induces isomorphisms on homotopy groups, and hence by Whitehead's 
theorem it is a homotopy equivalence.
Then the inclusion $j\:A\to MC(q)$ is also a homotopy equivalence.
Since $j$ is a cofibration, $MC(q)$ deformation retracts onto $A$.%
\footnote{Given a map $f\:MC(q)\to A$ and homotopies $h_t\:fj\simeq\id_A$ 
and $H_t\:jf\simeq\id_{MC(q)}$, let $r_t\:MC(q)\to A$ be an extension of $h_t$ over $MC(q)$ 
starting with $r_0=f$.
Then $r_1$ is a retraction of $MC(q)$ onto $A$ which is homotopic via $r_t$ and $H_t$
to $\id_{MC(q)}$.}
Therefore $MC(\Sigma^2p)$ equivariantly deformation retracts onto $\Sigma^2(X^2\x Y^\urcorner)$.
Hence $MC(\Sigma^2g)$ equivariantly deformation retracts onto $\Sigma^2f^{-1}(Y^\urcorner)$.
Then $MC(\Sigma^2g)\cup\Sigma^2Z$ equivariantly deformation retracts onto $\Sigma^2Z$.
On the other hand, the equivariant deformation retraction of $MC(f)$ onto $MC(g)\cup Z$ yields
an equivariant deformation retraction of $MC(\Sigma^2f)$ onto $MC(\Sigma^2g)\cup\Sigma^2Z$.
Thus $MC(\Sigma^2f)$ equivariantly deformation retracts onto $\Sigma^2Z$.
Consequently $\Sigma^2f\:\Sigma^2Z\to\Sigma^2(\tilde X*Y^\urcorner)$ is an equivariant 
homotopy equivalence.
If $\psi$ is its equivariant homotopy inverse, we obtain an equivariant map
$\Sigma^{k+2}\tilde X\xr{\Sigma^2\phi}\Sigma^2(Y^\urcorner*\tilde X)\xr{\psi}\Sigma^2Z
\xr{\chi}S^{m+k+1}$, where $\chi$ is the restriction of the given equivariant map 
$\widetilde{X\x CY}\to S^{m+k-1}$.
Since $2\dim X\le 2m-3$, by Lemma \ref{basic-b} it desuspends to an equivariant map 
$\tilde X\to S^{m-1}$.
\end{proof}

\subsection*{Acknowledgements}

I am very grateful to S. Parsa, R. Karasev and the two TMNA referees for useful discussions and feedback.
The geometric proof of Theorem \ref{t1} appeared in response to Parsa's talk at
the Moscow Geometric Topology Seminar.

\subsection*{Disclaimer}

I oppose all wars, including those wars that are initiated by governments at the time when they directly or indirectly support my research.
The latter type of wars include all wars waged by the Russian state since the Second Chechen war till the present day (April 22, 2022), as well as 
the US-led invasions of Afghanistan and Iraq.


\end{document}